\documentclass[11pt]{article}
\usepackage{amsmath, amssymb, chicago, amsfonts, latexsym, annals, natbib}
\advance \topmargin by -\headheight
\advance \topmargin by -\headsep
\advance \topmargin .2in
\def\E{\mathrm{E}}

\def\sz{S^{(0)}}

\def\bd{\begin{description}}
\def\ed{\end{description}}
\def\bl{\begin{list}{$\bullet$}{}}   % Begins a bullet list.
\def\cl{\begin{list}{$\circ$}{}}     % Begins a circle list.
\def\el{\end{list}}                  % Ends a list.

\def\b1{\mathbf{1}}                

\newcommand{\Var}{\mbox{Var}}
\newcommand{\PP}{{\mathbb  P}}
\newcommand{\QQ}{{\mathbb Q}}
\newcommand{\GG}{{\mathbb G}}
\newcommand{\RR}{{\mathbb R}}

\begin{document}
\setlength{\baselineskip}{20pt}

\title{Weighted Likelihood for Semiparametric Models
and Two-phase Stratified Samples, with Application to Cox Regression}

\author{Norman E. Breslow \\ Jon A. Wellner}

\affiliation{University of Washington, Seattle}
\date{\today}

\maketitle

\begin{abstract}
Weighted likelihood, in which one solves Horvitz-Thompson or inverse probability weighted (IPW) versions of the likelihood equations, offers a simple and robust method for fitting models to two phase stratified samples.
We consider semiparametric models for which solution of infinite dimensional estimating equations leads to $\sqrt{N}$ consistent and asymptotically Gaussian estimators of both Euclidean and nonparametric parameters.
If  the phase two sample is selected via Bernoulli (i.i.d.) sampling with known sampling probabilities, standard estimating equation theory shows that the influence function for the weighted likelihood estimator of the Euclidean parameter is the IPW version of the ordinary influence function.
By proving weak convergence of the IPW empirical process, and borrowing results on weighted bootstrap empirical processes, we derive a parallel asymptotic expansion for finite population stratified sampling.
Whereas the asymptotic variance for Bernoulli sampling involves the within strata second moments of the influence function, for finite population stratified sampling it involves only the within strata variances.
The latter asymptotic variance also arises when the observed sampling fractions are used as estimates of those known \textit{a priori}.
A general procedure is proposed for fitting semiparametric models with estimated weights to two phase data.
Several of our key results have already been derived for the special case of Cox regression with stratified case-cohort studies, other complex survey designs and missing data problems more generally.
This paper is intended to help place this previous work in appropriate context and to pave the way for applications to other models. 

\end{abstract}

\vspace{2cm}

\newpage
\section{Introduction}

Two phase stratified sampling, also known as double sampling, 
was introduced by \citet{4923} to estimate the population mean of a 
target variable that is costly or difficult to measure.
At phase one a relatively large random sample is drawn and 
measurements are made on an auxiliary variable that is 
correlated with the target variable but easier to measure.
At phase two measurements on the target variable are made
 for a subsample drawn randomly, without replacement, from 
 within strata defined by the auxiliary variable.
Neyman showed that the optimal, design unbiased linear 
estimator of the population mean is the Horvitz-Thompson (\citeyear{3888}) estimator that weights each observation by the inverse 
of the probability of its selection into the phase two sample.

Two-phase stratified sampling designs can dramatically reduce 
the costs of regression modeling when the strata depend 
on (correlates of) both outcome and explanatory variables.
A common method of estimation is ``weighted exogenous 
sampling maximum likelihood", here simply Weighted Likelihood 
or WL, in which one maximizes the inverse probability weighted 
(IPW) sum of log-likelihood contributions from the phase 
two observations \citep{200, 265}.
Equivalently, one may solve an IPW version of the 
score equations \citep[\S 3.4]{4105}.
Although easy to implement, WL estimators are sometimes 
seriously inefficient \citep{3937}.
They may still be of interest, however, because even when 
the model is wrong they consistently estimate the finite population 
parameters that would be obtained by fitting the model to 
complete phase one data \citep{1815, 4326}.
Fully efficient estimators are available for logistic and other 
parametric regression models in situations where the phase 
one data consist only of stratum frequencies.
See, for example, \citet{4927} and the references cited therein.

The asymptotic properties of WL estimators of Euclidean parameters 
in parametric models follow readily from standard results for 
$M$-estimators \citep[Chapter 5]{4895}.
WL may also be used for estimation of both Euclidean and infinite 
dimensional parameters in semiparametric models, for which the 
paradigm is Cox (\citeyear{1272}) proportional  hazards regression.  
\citet{4894} developed asymptotic results for both regression 
coefficients and baseline cumulative hazard when fitting the Cox model to survey 
data including those obtained using two phase sampling.
\citet{4324} obtained the same results for the regression parameters 
when fitting the Cox model to data from exposure stratified case-cohort 
studies, in which all subjects who have a failure event (the cases) 
are sampled at phase two.
One purpose of the present paper is to develop a modern theory
of WL estimation in semiparametric models that encompasses these previous results, helps to interpret them and paves the way for further applications.
We also explore the relationship between results based on finite 
population stratified sampling at phase two and those based on i.i.d. 
variable probability sampling with sampling weights 
estimated using information from phase one.

\section{Notation, Assumptions and Problem Statement}
Suppose $P_{\theta,\eta}$ denotes a probability distribution in 
a semiparametric model for a random variable $X \in {\cal X}$, 
where $\theta \in \Theta \subset \RR^p$ is the Euclidean parameter 
and $\eta$, taking values in some arbitrary space $H$, is the nonparametric one.
Let $P_0=P_{\theta_0,\eta_0}$ denote the distribution from which $X$ is actually sampled.
Following closely \S 25.12 of \citet{4895}, suppose maximum 
likelihood (ML) estimators $(\hat\theta, \hat\eta)$ are obtained by solving the system 
\begin{eqnarray}
\Psi_{N1}(\theta,\eta) & = & \PP_N \dot\ell_{\theta,\eta}\; = \;0 \nonumber \\
\Psi_{N2}(\theta,\eta) & = & \PP_N B_{\theta,\eta}h
      -P_{\theta,\eta}B_{\theta,\eta}h \;= \;0 \; \forall \;  h \in {\cal H}. \label{eq:like}
\end{eqnarray}
Here $\dot\ell_{\theta,\eta}$ is the $p$-dimensional 
likelihood score for $\theta$, $B_{\theta,\eta}$ is the score 
operator \citep{1262} working on an infinite dimensional class 
${\cal H}$ of directions $h$ from which paths of one-dimensional 
submodels for $\eta$ may approach $\eta_0$, and $\PP_N$ is empirical measure based on the i.i.d. sequence $X_1,\ldots,X_N$. 
Set $\dot\ell_0=\dot\ell_{\theta_0,\eta_0}$ and $B_0=B_{\theta_0,\eta_0}$.

Suppose the following assumptions, which slightly strengthen the 
hypotheses of \citet[Theorem 25.90]{4895}, are satisfied so that 
$\sqrt{N}(\hat\theta-\theta_0,\hat\eta-\eta_0)$ is asymptotically Gaussian:
\bd
\item[A1] for $(\theta,\eta)$ in a $\delta$-neighborhood of 
$(\theta_0,\eta_0)$   the functions $\dot\ell_{\theta,\eta}$ and 
 $\{B_{\theta,\eta}h, h \in {\cal H} \}$ are contained in a $P_0$-Donsker class ${\cal F}$;
\item[A2] $P_0\| \dot\ell_{\theta, \eta}-\dot\ell_0\|^2$ 
and $\sup_{h \in {\cal H}}P_0|B_{\theta, \eta }h-B_0h|^2$  converge  
to $0$ as $(\theta,\eta) \rightarrow (\theta_0,\eta_0)$;
\item[A3]
the map $\Psi=(\Psi_1,\Psi_2):\Theta \times H \mapsto \RR^p\times \ell^{\infty}({\cal H})$ 
with components
\begin{eqnarray}
\Psi_1(\theta,\eta) & = & P_0\dot\ell_{\theta,\eta} \nonumber \\
\Psi_2(\theta,\eta) & = & P_0B_{\theta,\eta}h
              -P_{\theta,\eta}B_{\theta,\eta} h , \ \  h \in {\cal H}, \label{eq:expectedmap}
\end{eqnarray}
which is the expectation of the random map $\Psi_N=(\Psi_{N1},\Psi_{N2})$ in (\ref{eq:like}), 
has a Fr\'echet derivative 
$\dot\Psi_0$ at $(\theta_0,\eta_0)$ that is continuously invertible on its range.
\item[A4] $(\hat\theta,\hat\eta)$ is consistent for $(\theta_0,\eta_0)$ and satisfies 
$\Psi_N(\hat\theta,\hat\eta)=0.$  
\ed
Assumption \textbf{A3} is typically established by showing that the information 
operator $B^*_0B_0$ is continuously invertible and thus that $\eta$ is 
estimable at a $\sqrt{N}$ rate.
This is the most restrictive assumption, but one that leads quickly to our main result.

With two phase sampling, however, $X$ is not observed for all $N$ subjects.
At phase one we observe only  a coarsening $\tilde X=\tilde X(X)$ of $X$ 
plus auxiliary variables $U \in {\cal U}$ that serve to determine the sampling strata.
$X$ is fully observed for subjects sampled at phase two.
Let $W=(X,U) \in {\cal W} = {\cal X} \times {\cal U}$ denote the variables 
potentially available for everyone, but in fact fully observed only for those 
in the phase two sample, and $V=(\tilde X,U) \in {\cal V} = {\cal \tilde X} \times {\cal U}$ 
denote the variables actually observed for everyone.  
We write $\tilde{P}_0$ for the distribution of $W = (X,U)$ and
denote by $\Sigma_N=\sigma[W_1,\ldots,W_N]$ the sigma field of information, 
also referred to as the complete data, potentially available for the $N$ subjects.
A sequence of binary indicators $(\xi_1,\ldots,\xi_N)$ shows which 
subjects are selected $(\xi_i=1)$ at phase two for observation of $X_i$.
We consider two probability models for the indicators $\xi_i$.
In the first, known as Bernoulli or Manski-Lerman (\citeyear{200}) sampling, 
each phase one subject is examined in succession for the value of 
$V_i$ and the indicator $\xi_i$ is independently generated with 
$\Pr(\xi_i=1|W_i)=\Pr(\xi_i=1|V_i) = \pi_0(V_i)$ where $\pi_0$ is a 
known sampling function. 
This preserves the i.i.d. structure for the observations $(\xi_i,V_i,\xi_iX_i)$.
Note the crucial missing at random (MAR) assumption: $\pi_0$ 
depends only on what is observed at phase one.  
We write $Q_0$ for the 
distribution of $(W_i, \xi_i)$.    
If ${\cal V}$ is partitioned into $J$ strata ${\cal V}_1 \cup \cdots \cup {\cal V}_J$, 
stratified Bernoulli sampling corresponds to the special case where 
$\pi_0(v)=p_j$ for $v \in {\cal V}_j$.
We assume that all $J$ strata are sampled with positive probability, or more generally that 
\begin{equation}
0 < \sigma \leq \pi_0(v) \leq 1 \quad \mbox{for} \quad v \in {\cal V}. 
\label{eq:boundedweights}
\end{equation}
Even though the sampling fractions are known, it is advisable to estimate 
$\pi_0$ in order to increase the efficiency of WL \citep{3937}.
We consider estimation of $\pi_0$ using a parametric model in \S 6.

The second sampling model corresponds to Neyman's original design 
and is usually closer to actual practice.
Here we observe the entire phase one sample at once and record the 
stratum frequencies 
$N_j=\sum_{i=1}^N \b1_{{\cal V}_j}(V_i)$ for $j=1,\ldots,J$.
At phase two samples of size $n_j \leq N_j$ are drawn at random, without 
replacement, from each of the $J$ finite phase one strata.
Using now a doubly subscripted notation where $\xi_{j,i}$ 
denotes the indicator variable for $i^{\mbox{th}}$ subject in stratum 
$j$, the essential features of this design are that, conditionally 
on $\Sigma_N$: ($i$) for $j=1,\ldots,J$ 
the random variables $(\xi_{j1},\ldots,\xi_{jN_j})$ are exchangeable 
with $\Pr(\xi_{j,i}=1|\Sigma_N)={n_j}/{N_j}$; 
and ($ii$) the $J$ random vectors $(\xi_{j1},\ldots,\xi_{jN_j})$ are independent.
Our problem is to estimate $(\theta,\eta)$ using the incomplete observations $V_i$ on everyone and the complete observations 
$X_i$ on subjects sampled at phase two.

\section{Weighted Likelihood Estimator}

WL estimates are obtained by solving Horvitz-Thompson 
(IPW) versions of the likelihood equations.
Define the \textit{inverse probability weighted empirical measure} by
\begin{equation}
\PP_N^{\pi} = \frac{1}{N}\sum_{i=1}^N\frac{\xi_i}{\pi_i}\delta_{X_i} , \label{eq:empiricalmeasure}
\end{equation}
where $\delta_{X_i}$ denotes Dirac measure placing unit mass on $X_i$ and
\begin{eqnarray*}
\pi_i & = & \left \{ \begin{array}{ll} \pi_0(V_i) & 
\mbox{\ for\ Bernoulli\ sampling}\\ \mbox{ } \\ \frac{n_j}{N_j} 
              \mbox{\ if\ }V_i\in{\cal V}_j &
\mbox{\ for\ finite\ population\ stratified\ sampling}. \end{array} \right .
\end{eqnarray*}
Then, instead of (\ref{eq:like}) we solve 
\begin{eqnarray}
\Psi_{N1}^{\pi}(\theta,\eta) & = & \PP_N^{\pi}\dot\ell_{\theta,\eta} \, = \,  0 \nonumber \\ 
\Psi_{N2}^{\pi}(\theta,\eta) & = & \PP_N^{\pi}B_{\theta,\eta}h - P_{\theta, \eta} B_{\theta, \eta} h 
\,  = \, 0  \qquad \mbox{for all} \ \ h \in {\cal H}. \label{eq:IPWlike}
\end{eqnarray}

In view of the MAR assumption, for any integrable function $f:{\cal X} \mapsto \RR$ 
and under either Bernoulli or finite population stratified sampling, 
\[
\E\frac{\xi_i}{\pi_i}f(X_i) = \E \left [ 
\E \left ( \left . \frac{\xi_i}{\pi_i} \right | \Sigma_N \right ) f(X_i) \right ] = \E f(X_i), \;  \; i=1,\ldots,N,
\]
so that $\E \PP_N^{\pi}f = \E \PP_N f = P_0f$.
Consequently, the random map $\Psi_N^{\pi}=(\Psi_{N1}^{\pi},\Psi_{N2}^{\pi})$ 
defined by (\ref{eq:IPWlike}) has the same expectation as the random 
map $\Psi_N$ in (\ref{eq:like}), namely $\Psi=(\Psi_1,\Psi_2)$ as in (\ref{eq:expectedmap}).
The implication is that the assumptions \textbf{A1}-\textbf{A4} 
made to guarantee the asymptotic normality of the ML estimator 
based on complete phase one data are also the assumptions needed 
to guarantee the asymptotic normality of the WL estimator based on two phase data.
Indeed, van der Vaart's (\citeyear{4895}) Theorem 25.90, 
or more precisely his Theorem 19.26 of which it is a restatement, 
applies virtually without change to the Bernoulli sampling setup. 
The Donsker class ${\cal F}$ in \textbf{A1} is modified to
$\tilde{\cal F} = \{ [\xi/\pi_0(V)] f(X), f \in {\cal F}\}$. 
Since under the hypothesis (\ref{eq:boundedweights}) it is the product 
of a fixed bounded function with the Donsker class ${\cal F}$, the fact
that $\tilde{\cal F}$ is Donsker for the joint distribution $Q_0$ 
of $(W, \xi)$ follows from \citet[example 2.10.10]{4920}.
The random map $\Psi_N$ corresponding to the estimating 
functions (\ref{eq:IPWlike}) is ordinary empirical measure 
$\QQ_N$ for $\{(W_i, \xi_i), i=1,\ldots, N\}$ applied to the unbiased 
estimating functions $(\xi/\pi_0)\dot\ell_{\theta,\eta}$ and $(\xi/\pi_0)B_{\theta,\eta}h$.
\textbf{A4} will generally follow from (\ref{eq:boundedweights}) 
and the arguments used to establish consistency for the complete data ML estimator, 
together with (\ref{eq:IPWlike}).  
\textbf{A2} and \textbf{A3} are unchanged.
The more general Theorem 3.3.1 of \citet{4920} is needed, however, 
to deal with the non i.i.d. data induced by finite population stratified sampling.
To verify its hypotheses, we first must establish weak convergence 
of the empirical process based on $\PP_N^{\pi}$.

\section{Weak Convergence of the IPW Empirical Process}

Two phase stratified sampling resembles the bootstrap in that it 
involves random sampling from the finite, albeit incompletely 
observed, population $\{X_1,\ldots,X_N\}$.
Here we use results on weighted bootstrap empirical processes 
from \citet[Theorem 2.2]{4922}, as incorporated in \citet[Theorem 3.6.13]{4920}, 
to demonstrate weak convergence of the IPW empirical process 
$\mathbb{G}_N^{\pi}=\sqrt{N}(\PP_N^{\pi}-P_0)$ for finite population stratified sampling.
First note that, with the subscript $j,i$ denoting the 
$i^{\mbox{th}}$ of $N_j$ observations in stratum $j$, 
\begin{eqnarray}
\mathbb{P}_N^{\pi} & = & \frac{1}{N}\sum_{j=1}^J\frac{N_j}{n_j}\sum_{i=1}^{N_j}\xi_{j,i}\delta_{X_{j,i}} 
 =  \frac{1}{N}\sum_{j=1}^J \frac{N_j^2}{n_j} \mathbb{P}_{j,N_j}^{\xi} \label{eq:bootstrap}
\end{eqnarray}
where 
\[
\mathbb{P}_{j,N_j}^{\xi} = \frac{1}{N_j}\sum_{i=1}^{N_j}\xi_{j,i}\delta_{X_{j,i}} 
\]
is a \textit{finite sampling empirical measure} for the $j^{\mbox{th}}$ 
stratum.
Similarly one can express the ordinary empirical measure as
\begin{equation}
\mathbb{P}_N = \frac{1}{N}\sum_{j=1}^JN_j\mathbb{P}_{j,N_j} 
\label{eq:empirical}
\end{equation}
where 
\begin{equation}
\mathbb{P}_{j,N_j} = \frac{1}{N_j} \sum_{i=1}^N \delta_{X_i} \b1_{{\cal V}_j} (V_i) 
= \frac{1}{N_j}\sum_{i=1}^{N_j}\delta_{X_{j,i}}
\label{eq:stratumwiseEmpirical}
\end{equation}
denotes the empirical measure for the $j^{\mbox{th}}$ stratum.  
Justification of the second (doubly indexed) form is given in Appendix A.

Combining (\ref{eq:bootstrap}) and (\ref{eq:empirical}), and letting 
$\GG_N=\sqrt{N}(\PP_N-P_0)$ denote the standard empirical process, we have
\begin{eqnarray}
\mathbb{G}_N^{\pi} 
& = & \sqrt{N}\left(\mathbb{P}_N^{\pi}-P_0\right) \nonumber \\
& = & \sqrt{N}\left(\mathbb{P}_N-P_0\right)
           +\sqrt{N}\left(\mathbb{P}_N^{\pi}-\mathbb{P}_N\right) \nonumber \\
& = & \GG_N + \frac{1}{\sqrt{N}}\sum_{j=1}^J
             \left(\frac{N_j^2}{n_j}\right)\left(\mathbb{P}_{j,N_j}^{\xi}-\frac{n_j}{N_j}
             \mathbb{P}_{j,N_j}\right) \nonumber \\
& = & \GG_N 
            +\sum_{j=1}^J\sqrt{\frac{N_j}{N}}\left(\frac{N_j}{n_j}\right)\mathbb{G}^{\xi}_{j,N_J} 
            \label{eq:IPWempiricalexpansion}
\end{eqnarray}
where
\begin{equation}
\mathbb{G}^{\xi}_{j,N_j} = \sqrt{N_j}\left(\mathbb{P}_{j,N_j}^{\xi}
       -\frac{n_j}{N_j}\mathbb{P}_{j,N_j}\right) \label{eq:weightedbootstrapprocess}
\end{equation}
is the \textit{finite sampling empirical process} for stratum $j$.

The first term in (\ref{eq:IPWempiricalexpansion}) converges to the 
$P_0$-Brownian bridge process $\mathbb{G}$ indexed by the Donsker 
class ${\cal F}$ mentioned in \textbf{A1}.
Let $P_{0|j}(\cdot)=\E(\cdot|V \in{\cal V}_j)$ denote $\tilde{P}_0$ 
conditional on membership in stratum $j$, \textit{i.e.},  for measurable 
$A \subset {\cal X}$, $P_{0|j} (A) = \tilde P_0 [A  {\bf 1}_{{\cal V}_j}(V)]/\nu_j$ 
with $\nu_j= \tilde P_0\b1_{{\cal V}_j}(V)$,
and let $\mathbb{G}_j$ denote the 
$P_{0|j}$-Brownian bridge, also indexed by $ {\cal F}$.
Our goal is to establish the weak convergence of the remaining terms 
on the RHS of (\ref{eq:IPWempiricalexpansion}).
If as $N \rightarrow \infty$ the sampling fractions converge 
with $n_j/N_j \rightarrow p_j$,  the assumption on the exchangeable 
``weights" $(\xi_{j,1},\ldots,\xi_{j,N_j})$ in equation (3.6.8) of \citet{4920} holds trivially with
\[
\frac{1}{N_j}\sum_{i=1}^{N_j}\left(\xi_{j,i}-\bar\xi_{j.}\right)
\stackrel{\mbox{p}}{\rightarrow} p_j(1-p_j) \label{eq:VdVW368}.
\]
Furthermore, with $\rightsquigarrow$ denoting weak convergence in 
$\ell^{\infty}({\cal F})$, $\sqrt{N_j} ( \PP_{j,N_j} - P_{0|j} ) \rightsquigarrow \GG_{j}$; 
see Appendix B for the proof.
Thus their Theorems 3.6.13 and 1.12.4 imply that, for almost 
every sequence of complete data, 
$\mathbb{G}_{j,N_j}^{\xi} \rightsquigarrow \sqrt{p_j(1-p_j)}\mathbb{G}_j$.
Conditionally on $\Sigma_N$, the processes $\mathbb{G}^{\xi}_{j,N_j}$ 
are mutually independent because of the independence of the 
$\{\xi_{j,i}\}$ in different strata. 
Furthermore, by virtue of the fact that they also are (unconditionally) 
uncorrelated with $\mathbb{G}_N=\sqrt{N}(\mathbb{P}_N-P_0)$, 
which follows along the lines of \citet[Corollary 2.9.3]{4920}, 
or that (conditionally) they have the same limiting distributions for 
almost all sequences of data, the vector of processes 
$(\mathbb{G}_N,\mathbb{G}_{1,N_1}^{\xi},\ldots,\mathbb{G}_{J,N_J}^{\xi})$ 
converges weakly to the vector of independent Brownian bridge processes
$(\mathbb{G},\mathbb{G}_1,\ldots,\mathbb{G}_J)$.
Consequently
\begin{equation}
\mathbb{G}_N^{\pi} \rightsquigarrow \mathbb{G} 
+ \sum_{j=1}^J\sqrt{\nu_j}\sqrt{\frac{1-p_j}{p_j}} \mathbb{G}_j. 
\label{eq:limIPWempiricalprocess}
\end{equation}
This result formalizes and extends Proposition 1 of \citet{218} 
and the arguments in \S 4 of \citet{4324}.

\section{Asymptotic Distributions of the WL estimator}

We apply Theorem 19.26 of \citet{4895} to conclude that, 
under Bernoulli sampling, 
\begin{equation}
\sqrt{N} \dot\Psi_0 \left( \begin{array}{c}
\hat{\theta}-\theta_0 \\ \hat{\eta}-\eta_0 \end{array} \right ) 
= -\GG_N \frac{\xi}{\pi_0}
\left (  \begin{array}{c} \dot\ell_0 \\ 
 B_0h \end{array} \right ) \; + \; o_p(1). \label{eq:BernoulliForm}
\end{equation}
Similarly, using Theorem 3.3.1 of \citet{4920} together with the development 
of the previous section, we conclude that 
for finite population stratified sampling 
\begin{equation}
\sqrt{N} \dot\Psi_0 \left( \begin{array}{c}
\hat{\theta}-\theta_0 \\ \hat{\eta}-\eta_0 \end{array} \right ) 
= -\GG_N^{\pi} 
\left ( \begin{array}{c} \dot\ell_0 \\  
B_0h \end{array} \right ) \;+ \;o_p(1). \label{eq:FPSSForm}
\end{equation}
We have already argued that the hypotheses of the first theorem 
follow from appropriately modified versions of \textbf{A1}-\textbf{A4}.
Together with the weak convergence of $\mathbb{G}_N^{\pi}$ 
just established, they also suffice for the second theorem.
In particular, the stochastic condition (3.3.2) of \citet{4920} follows 
from \textbf{A1} and \textbf{A2} together with the proof of their Lemma 3.3.5 
applied to each of $\GG_N,\GG^{\xi}_{1,N_1},\ldots, \GG^{\xi}_{1,N_1}$.

In practice attention is usually focused on inferences for the Euclidean parameter $\theta$. 
To derive a general expression for the asymptotic variance of $\hat\theta$ we further assume 
\bd
\item[A5] $\dot\Psi_0$ admits a partition as in equation (25.91) 
of \citet{4895} where the information operator 
$B_0^* B_0$ is continuously invertible.
\ed
Following closely the arguments in \S 25.12 of van der Vaart, we 
calculate from (\ref{eq:BernoulliForm}) that under Bernoulli sampling 
\begin{equation}
\sqrt{N}(\hat{\theta}-\theta_0) = \GG_N \frac{\xi}{\pi_0}\tilde\ell_0 + o_p (1) \label{eq:mainresultBernoulli}
\end{equation}
whereas from (\ref{eq:FPSSForm}) under finite population stratified sampling
\begin{equation}
\sqrt{N}(\hat{\theta}-\theta_0)= \GG_N^{\pi} \tilde\ell_0 + o_p (1),
\label{eq:mainresultFPSS}
\end{equation}
where in both cases $\tilde\ell_0$ denotes the efficient influence function
\begin{equation}
\tilde\ell_0 = \tilde I_0^{-1}\left ( I-B_0\left ( B^*_0B_0 \right )^{-1}  B^*_0 \right ) \dot\ell_0 \label{eq:effinfluence}
\end{equation}
and 
\begin{equation}
\tilde I_0 
= P_0 \left [  \left ( I-B_0\left ( B^*_0B_0 \right )^{-1}  B^*_0 \right ) 
\dot\ell_0\dot\ell_0^{T} \right ] \label{eq:effinfo}
\end{equation}
is the efficient information.
Since $P_0\tilde\ell_0=0$, moreover, both (\ref{eq:mainresultBernoulli}) and (\ref{eq:mainresultFPSS}) may be expressed 
\begin{equation}
\sqrt{N}(\hat\theta-\theta_0) = \sqrt{N} \PP_N^{\pi}\tilde\ell_0 + o_p(1)= 
\frac{1}{\sqrt{N}} \sum_{i=1}^N \frac{\xi_i}{\pi_i}\tilde\ell_0(X_i) +o_p(1), \label{eq:general}
\end{equation}
which expansion constitutes the principal result of this paper.

Under Bernoulli sampling with known $\pi_0$ the asymptotic variance is therefore
\begin{eqnarray}
\Var_{\mbox{A}}\sqrt{N}(\hat\theta-\theta_0) & = & 
\Var \left ( \frac{\xi}{\pi_0} \tilde\ell_0 \right ) \nonumber \\
 & = & 
\Var \; \E \left ( \left . \frac{\xi}{\pi_0} \tilde\ell_0 \right | X \right ) + 
\E \; \Var \left ( \left . \frac{\xi}{\pi_0} \tilde\ell_0 \right | X \right ) \nonumber \\
& = & \Var(\tilde\ell_0) + \E \left [ \frac{\tilde\ell_0^{\otimes 2}}{\pi_0^2} \Var(\xi|X)\right ]
\nonumber \\
& = & \tilde I_0^{-1} + \tilde P_0 \left ( \frac{1-\pi_0}{\pi_0}\tilde\ell_0^{\otimes 2} \right ). \label{eq:IPWvariance}
\end{eqnarray}
In the special case of stratified Bernoulli sampling, with $\pi_i=\pi_0(V_i)=p_j$ for $V_i \in {\cal V}_j$, this becomes
\begin{equation}
\tilde I_0^{-1} + \sum_{j=1}^J\nu_j\frac{1-p_j}{p_j}P_{0|j}\left (\tilde\ell_0^{\otimes 2} \right ).
\label{eq:varIPWstrataiid}
\end{equation}
On the other hand, from (\ref{eq:limIPWempiricalprocess}) and 
(\ref{eq:mainresultFPSS}), the asymptotic variance under finite population stratified sampling is 
\begin{equation}
\tilde I_0^{-1} + \sum_{j=1}^J\nu_j\frac{1-p_j}{p_j}\Var_j(\tilde\ell_0), \label{eq:varIPWstratified}
\end{equation}
where $\Var_j(f)=P_{0|j}(f^{\otimes 2})-P^{\otimes 2}_{0|j}(f)$.
Comparing the last two expressions shows the substantial potential gain from 
keeping track of the stratum frequencies for the phase one data.

\section{Bernoulli Sampling with Estimated Weights}
Let ${\cal V}_0$ denote an additional stratum, possibly null, such that $\xi_i=1$ for $V_i \in {\cal V}_0$.
Introduction of this special stratum with $p_0=1$ does not affect the previous development; 
in particular, equations (\ref{eq:general})-(\ref{eq:varIPWstratified}) continue to hold.
For $V_i \notin {\cal V}_0$ suppose 
\begin{equation}
\Pr(\xi_i=1|X_i, V_i;\alpha) = \Pr(\xi_i=1|V_i;\alpha) = \pi_{\alpha}(V_i)  < 1 \label{eq:modelforPi}
\end{equation}
where $\alpha \in \Xi \subset \RR^q$ 
is a parameter to be estimated by 
maximum likelihood from the phase 
one observations $\{V_i, i=1,\ldots,N\}$ not in ${\cal V}_0$.
We assume sufficient regularity in the model for $\alpha$, e.g., to satisfy the hypotheses of 
Theorem 5.21 of \citet{4895}, so that the ML estimator 
$\hat\alpha$ is consistent and asymptotically normal with influence function
\begin{equation}
\tilde\ell^{\alpha}_0= \b1_{{\cal V}_0^c} 
\left ( \tilde P_0 \b1_{{\cal V}_0^c} \frac{\dot\pi_0^{\otimes 2}} {\pi_0(1-\pi_0)} \right )^{-1} \dot\pi_0\frac{\xi-\pi_0}{\pi_0(1-\pi_0)}. 
\label{eq:InfFuncAlpha}
\end{equation}
Here for $V \in {\cal V}_0^c$, the complement of ${\cal V}_0$,  
$\pi_0(V)=\pi_{\alpha_0}(V)$ is the true sampling function while $\dot\pi_0(V)$ 
denotes the $q$-vector of partial derivatives of $\pi_{\alpha}(V)$ with 
respect to $\alpha$ evaluated at $\alpha=\alpha_0$.
If $\hat\theta(\alpha)$ denotes the WL estimator under two phase 
Bernoulli sampling with ``known" sampling function $\pi_{\alpha}(V)$, 
then from (\ref{eq:InfFuncAlpha}) and (\ref{eq:general}) we have
\begin{equation}
\sqrt{N}\left(\begin{array}{c}\hat\theta(\alpha_0)-\theta_0\\[.2cm] 
\hat\alpha-\alpha_0 \end{array}\right) 
= \sqrt{N}\left(\begin{array}{c}\PP_N^{\pi}\tilde\ell_0 \\[.2cm] \QQ_N\tilde\ell^{\alpha}_0
\end{array}\right) + o_p(1). \label{eq:jointExpansion}
\end{equation}
Furthermore, with $\hat\pi_i=\pi(V_i;\hat\alpha)$ for 
$V_i \in {\cal V}_0^c$ otherwise $\hat\pi_i=1$, we show in Appendix C that under some further mild assumptions regarding $\pi_{\alpha}(V)$
%\textbf{[Jon: please note changes, in particular use 
%of $\QQ$ and reference to additional mild assumptions in Appendix]}
\begin{equation}
\sqrt{N}(\PP_N^{\hat\pi}-\PP_N^{\pi_0})\tilde\ell_0  
= - \tilde  P_0 \left (\b1_{{\cal V}_0^c} \frac{\tilde\ell_0 
\dot\pi^{T}_0}{\pi_0}\right ) \sqrt{N}(\hat\alpha-\alpha_0) + o_p(1). \label{eq:TaylorExpansion} 
\end{equation}
The joint asymptotic normality of $(\hat\theta(\alpha_0),\hat\alpha)$ that follows
from (\ref{eq:jointExpansion}), together with the Taylor expansion 
(\ref{eq:TaylorExpansion}), are precisely the hypotheses used by \citet{855} 
to deduce that $\sqrt{N}[\hat\theta(\hat\alpha)-\theta_0] \rightsquigarrow Z $ 
where  $Z\in \RR^p $ is mean zero Gaussian with covariance matrix
\begin{equation}
\Var_{\mbox{A}}\sqrt{N}\left(\hat\theta(\hat\alpha)-\theta_0\right) 
= \Var  \left ( \frac{\xi}{\pi_0}\tilde\ell_0 \right ) - \tilde P_0 \b1_{{\cal V}_0^c}
    \frac{\tilde\ell_0 \dot\pi_0^T}{\pi_0} \left ( \tilde P_0 \b1_{{\cal V}_0^c}
    \frac{\dot\pi_0^{\otimes 2}}{\pi_0(1-\pi_0)} \right )^{-1} 
    \tilde P_0 \b1_{{\cal V}_0^c} \frac{\dot\pi_0\tilde\ell_0^T }{\pi_0}. \label{eq:VarEstAlpha}
\end{equation}
A matrix calculation shows that, when (\ref{eq:VarEstAlpha}) is evaluated for stratified Bernoulli sampling 
\[
 \pi_{\alpha} = \pi_{\alpha}(V) = \left \{ \begin{array}{ll} 1,  & V \; \in \; {\cal V}_0 \\ 
 \alpha_j, & V \; \in \; {\cal V}_j, \; j=1,\ldots,J , \end{array} \right . 
\]
the asymptotic variance for the WL estimator $\hat\theta$ with 
\textit{estimated} sampling probabilities $\hat\alpha_j=n_j/N_j$ is identical to 
the finite population sampling variance (\ref{eq:varIPWstratified}) with $p_j=\alpha_{j,0} = \lim n_j/N_j$.

Two possibilities present themselves for estimation of the terms in (\ref{eq:VarEstAlpha}).
Let $\hat\pi_i= \pi_{\hat\alpha}(V_i)$ for $V_i \in {\cal V}_0^c$ else $\hat\pi_i=1$.
Then, using (\ref{eq:IPWvariance}), we could estimate the first term by
\[
\widehat{\Var\left ( \frac{\xi}{\pi_0}\tilde\ell_0 \right )} 
= \tilde I^{-1}_{\hat\theta,\hat\eta} + \frac{1}{N}\sum_{i=1}^N\frac{\xi_i(1-\hat\pi_i)}{\hat\pi_i^2}\tilde\ell_{\hat\theta,\hat\eta}^{\otimes 2}(X_i),
\]
the expression in the middle of the second term by
\[
\widehat {\tilde P_0 \b1_{{\cal V}_0^c} \frac{\dot\pi_0^{\otimes 2}}{\pi_0(1-\pi_0)} }
= \frac{1}{N}\sum_{i=1}^N \b1_{{\cal V}_0^c}(V_i) 
\frac{\dot\pi_{\hat\alpha}^{\otimes 2}(V_i)}{\hat\pi_i(1-\hat\pi_i)}
\]
and similarly for $\tilde P_0(\tilde\ell_0\dot\pi_0^{T}/\pi_0)$.
A more empirical approach, however, would be to use the 
$\theta$ and $\alpha$ influence function contributions themselves to estimate these terms as in
\begin{eqnarray*}
\widehat{\Var\left ( \frac{\xi}{\pi_0}\tilde\ell_0 \right ) }
& = & \frac{1}{N}\sum_{i=1}^N \left ( \frac{\xi_i}{\hat\pi_i}\tilde\ell_{\hat\theta,\hat\eta}(X_i) \right )^{\otimes 2}, \\
\widehat {\tilde P_0 \b1_{{\cal V}_0^c} \frac{\tilde\ell_0\dot\pi_0^{T} } {\pi_0} } 
& = & \frac{1}{N}\sum_{i=1}^N \b1_{{\cal V}_0^c}(V_i) \frac{\xi_i}{\hat\pi_i} 
          \frac{\tilde\ell_{\hat\theta,\hat\eta}(X_i)} {\hat\pi_i} \dot\pi_{\hat\alpha}(V_i)^{T} \\
& = & \frac{1}{N}\sum_{i=1}^N \b1_{{\cal V}_0^c }(V_i) 
           \left ( \frac{\xi_i\tilde\ell_{\hat\theta,\hat\eta}(X_i)}{\hat\pi_i} \right ) 
           \left ( \frac{\dot\pi_{\hat\alpha}(V_i)^{T}(\xi_i-\hat\pi_i)}{\hat\pi_i(1-\hat\pi_i)} \right ) \quad \mbox{and} \\
\widehat { \tilde P_0 \b1_{{\cal V}_0^c} \frac{\dot\pi_0^{\otimes 2}}{\pi_0(1-\pi_0)} } 
& = & \frac{1}{N}\sum_{i=1}^N \b1_{{\cal V}_0^c }(V_i) 
            \left ( \frac{ \dot\pi_{\hat\alpha}(V_i)(\xi_i-\hat\pi_i)}{\hat\pi_i(1-\hat\pi_i)} \right )^{\otimes 2}.
\end{eqnarray*}
The resulting asymptotic variance for $\hat\theta$ may be recognized as the comprising the residual sums of squares and of cross products from the least squares regressions of each the $p$ components of the $\hat\theta$ influence function contributions $\xi_i\tilde\ell_{\hat\theta,\hat\eta}(X_i)/\hat\pi_i$, to which subjects not in the 
phase two sample contribute 0, on the $q$ components of the estimated $\hat\alpha$ influence function contributions (\ref{eq:InfFuncAlpha}), to which subjects having $V_i \in {\cal V}_0$ contribute 0. 
See \citet{4931} for a recent discussion and interpretation.
This suggests the following estimation procedure:
\begin{enumerate}
\item Estimate $\alpha$ from the phase one data and compute the estimated sampling fractions $\hat\pi_i$.
\item Estimate $\theta$ and $\eta$ from the phase two data by WL, using the inverse $\hat\pi_i$ as known weights.
\item Regress each component of the influence function contributions for $\hat\theta$ on those for $\hat\alpha$.
\item Estimate Var$_{\mbox{A}}(\hat\theta)$ as the matrix comprising the residual sums of squares and of cross products from these regressions.
\end{enumerate}
\citet[p. 166]{4901}, who cited earlier work by \citet{4902}, suggested this 
procedure for the special case of Cox regression, to which we now direct our attention.

\section{Application to the Cox Proportional Hazards Model}
Our development of the Cox model follows closely that of \citet[\S 25.12]{4895} 
where $X=(\Delta,T,Z)$ with $T$=min($\tilde T,C)$ a censored failure time, 
$\Delta=\b1_{[\tilde T \leq C]}$ the failure indicator and $Z \in \RR^p$ a vector of covariates. 
The Euclidean parameter is the $p$-vector of regression coefficients $\beta$ in the linear predictor $z\beta$.
The nonparametric parameter $\eta=(\Lambda,G,G_Z)$ has three 
infinite dimensional components: $\Lambda(\cdot)=\int_0^{\cdot}\lambda(s)ds$ the baseline cumulative 
hazard function, assumed differentiable; $G(t|z)=\Pr(C \leq t|Z=z)$ the conditional distribution of the 
censoring time; and $G_Z$, the marginal distribution of the covariates.
We introduce the usual notation for the ``at risk" process 
$Y(t)=\b1_{[T \geq t]}$ and the event counting process $N(t)=\Delta \b1_{[T \leq t]}$ 
and we make the standard assumptions: (i) that the true failure time $\tilde T$ 
and $C$ are independent given $Z$; and (ii) that there is a finite maximum 
censoring time $\tau$ such that $\Pr[Y(\tau)=1]>0$.
\citet{4895} makes some further ``partly unnecessary" assumptions to simplify his development, namely that the covariates $Z$ are bounded,  
that $G$ and $G_Z$ have densities as indicated and especially that $\Pr(C \geq \tau) = \Pr(C = \tau) >0$ (see discussion in \S 8).
Writing the density for $x=(\delta,t,z)$, with $z$ a row vector, as
\begin{equation}
e^{-e^{z\beta}\Lambda(t)}\left [e^{z\beta}\lambda(t)\left (1-G(t-|z)\right)\right]^{\delta} 
\left[g(t|z)\right]^{1-\delta}g_Z(z) \label{eq:Coxlike},
\end{equation}
and noting that $G$ and $G_Z$ factor out of the complete data likelihood, \citet{4895} 
considers ML estimation for $(\beta,\Lambda)$ only.
With ${\cal H}$ denoting various subsets of the space BV$[0,\tau]$ of bounded functions of bounded variation,
he develops the following explicit expressions for the $\beta$ score vector, 
the $\Lambda$ score operator that maps functions $h \in {\cal H}$ to functions of the data, its adjoint (but only evaluated for the 
$\beta$ scores) and the information operator that maps ${\cal H}$ onto itself:
\begin{eqnarray}
\dot\ell_{\beta,\Lambda}(x) 
& = & \delta z -ze^{z\beta}\Lambda(t) \label{eq:Coxscore} \\
B_{\beta,\Lambda}h(x) 
& = & \delta h(t) - e^{z\beta}\int_0^thd\Lambda \label{eq:CoxScoreOp} \\
B^*_{\beta,\Lambda}\dot\ell_{\beta,\Lambda}(t) 
& = & P_{\beta,\Lambda}Y(t)Ze^{Z\beta} \nonumber \\
B^*_{\beta,\Lambda}B_{\beta,\Lambda}h(t) 
& = & h(t) P_{\beta,\Lambda}Y(t)e^{Z\beta}. \nonumber
\end{eqnarray}
These are used to calculate the efficient scores 
\begin{eqnarray*}
\ell^*_{\beta,\Lambda}(x)& = & \dot\ell_{\beta,\Lambda}-
B_{\beta,\Lambda}\left(B^*_{\beta,\Lambda}B_{\beta,\Lambda}\right)^{-1}B^*_{\beta,\Lambda}\dot\ell_{\beta,\Lambda} \\
 & = & \delta\left[z-m(t;\beta)\right]-e^{z\beta}\int_0^t\left[z-m(s;\beta)\right]d\Lambda(s) 
\end{eqnarray*}
and efficient information
\begin{eqnarray*}
 \tilde I_0 & = & I_0 - P_0B_0\left(B^*_0B_0\right)^{-1}B^*_0\dot\ell_0 \\
 & = & P_0 \left ( e^{Z\beta_0}\int_0^{\tau} \left [
Z - m(t;\beta_0) \right ]^{\otimes 2} \Pr(T \geq t |Z) d\Lambda_0(t) \right ) ,
\end{eqnarray*}
respectively, 
where $I_0=P_0\dot\ell_o\dot\ell_0^{T}$ and $m(t;\beta) = S^{(1)}(t;\beta)/S^{(0)}(t;\beta)$ with
\begin{eqnarray*}
S^{(0)}(t;\beta) &=& P_0 e^{Z\beta}Y(t)  \\
S^{(1)}(t;\beta) &=& P_0Ze^{Z\beta}Y(t).
\end{eqnarray*}

To fit the Cox model by WL to two phase stratified samples,
first define IPW estimators of the two quantities just considered by
$\hat S^{(0)}(t;\beta)=\PP_N^{\pi}e^{Z\beta}Y(t)$ and 
$S^{(1)}(t;\beta)=\PP_N^{\pi}Ze^{Z\beta}Y(t) $.
By definition the WL estimators solve
\begin{eqnarray}
\Psi_{N1}^{\pi}(\beta,\Lambda) & = & \PP_N^{\pi}\dot\ell_{\beta,\Lambda} = 0 
\label{eq:CoxWL1} \\
\Psi_{N2}^{\pi}(\beta,\Lambda)h & = & \PP_N^{\pi} B_{\beta,\Lambda}h =0
\qquad \mbox{for all} \; h\in \; {\cal H}, \label{eq:CoxWL2}
\end{eqnarray}
where we have used the fact that $P_{\beta,\Lambda}B_{\beta,\Lambda}h=0 $. 
Substituting 
\[ 
h_t(s) \; = \; \frac{\b1_{[s\leq t]}}{\hat S^{(0)}(s,\beta)} 
\]
for $h$ in (\ref{eq:CoxWL2}) and solving using (\ref{eq:CoxScoreOp}) 
shows that, for fixed $\beta$, the cumulative hazard 
function that partially maximizes the weighted likelihood and, as is easily checked, satisfies $\PP_N^{\pi}B_{\beta,\hat\Lambda_{\beta}}h=0$ for all  $h$, is
\begin{equation} 
\hat\Lambda_{\beta}(t) \;= \; \PP_N^{\pi}\frac{\Delta\b1[T\leq t]}{\hat\sz(T;\beta)} \; = \;\frac{1}{N} \sum_{i=1}^N \int_0^t\frac{\xi_i}{\pi_i}\frac{dN_i(s)}{\hat S^{(0)}(s;\beta)} \label{eq:Breslow}.
 \end{equation}
This may be recognized as an IPW version of the so called \citet{1266} estimator.
Inserting this expression into (\ref{eq:CoxWL1}) and evaluating using (\ref{eq:Coxscore}) yields
\[
\Psi_{N1}^{\pi}(\beta,\hat\Lambda_{\beta})\; = \; \PP_N^{\pi}\Delta\left[Z-\hat m(T;\beta)\right] \;= \; \frac{1}{N}\sum_{i=1}^N\frac{\xi_i}{\pi_i}\Delta_i\left[Z_i-\frac{\hat S^{(1)}(T_i;\beta}{\hat S^{(0)}(T_i;\beta)}\right]\;=\;0 ,
\]
which is the IPW Cox ``partial score" equation.
Its solution, together with (\ref{eq:Breslow}),
are the estimators proposed for Cox regression by \citet{4326}, 
\citet{4902}, \citet[Estimator II]{4324}, \citet{4894} and others for a 
variety of complex sampling and missing data problems.
Using the results of this paper, the large sample properties of $(\hat\beta,\hat\Lambda_{\hat\beta})$ follow 
from those already developed for the ML estimators with complete data, which are given by the same equations with $\xi_i=\pi_i=1, i=1,\ldots,N$.

\section{Discussion}
The two phase stratified sampling designs considered here are quite 
flexible in that the phase one strata may be formed using all available 
information and sampled with arbitrary positive probabilities.
This is in the spirit of \citet{4326} and \citet{4894}, who considered 
even more general complex sample survey designs.
Others \citep{4324,4871} have restricted their attention to covariate stratified versions of the case-cohort design, whereby all subjects who fail are sampled at phase two for complete covariate ascertainment.
Although this may well be an efficient design when the failure rate is low, the assumption that $\xi=1$ whenever $\Delta=1$ is often unnecessary and may sometimes be unduly restrictive.
Not only does it limit application when the phase one population has 
large numbers of both failures and non-failures, it also does so when the sampling has been carried out for one failure type but it is of interest to evaluate another.
When following patients enrolled in a clinical trial, for example, all deaths may be sampled as ``cases" but it may later be decided to analyze the data also in terms of ``event-free survival".
In other contexts, biological samples may turn out out to be non-informative so that data are still missing for substantial numbers of subjects, including failed cases, who are sampled at phase two.
Provided one is willing to make the standard MAR assumptions, WL 
methods as described herein may still be used by determining the 
stratum frequencies for subjects having complete data at phase two and using them to estimate the sampling weights.

The major drawback of WL estimation is its lack of statistical efficiency.
Efforts to address this deficiency with Cox regression have been made by 
several authors including \citet{3937}, \citet{4871}, \citet{4811}, \citet{4929} and \citet{4930}.
Most of these methods are relatively recent and involve sufficiently complex calculations, or sufficiently restrictive assumptions, that none have yet seen widespread use.
These limitations are certain to decline with advances in computing 
hardware and software, making more efficient estimation methods more widely available.
In the meantime, the WL estimation procedure outlined at the end 
of \S 6 offers a relatively simple and robust alternative.
It is likely to remain the method of choice for many survey statisticians for the reasons 
mentioned in the introduction, namely, their interest in finite population 
parameters defined as solutions to ML estimating equations.
As emphasized by \citet{3937}, in view of the interpretation of 
(\ref{eq:VarEstAlpha}) as a residual sum of squares, inclusion of 
additional variables in the model (\ref{eq:modelforPi}) for $\pi$ can 
only enhance the efficiency of $\theta$ estimation.
When the sampling probabilities vary, as in finite population stratified 
sampling, inclusion of the stratum factors in the model is essential to avoid bias.
Finer stratification, or the inclusion of auxiliary variables in the model for $\pi$, serves the cause of efficiency.
Equation (\ref{eq:varIPWstratified}) suggests that such additional 
variables would be most valuable if they could somehow be chosen 
to be highly correlated with the efficient scores. 
The doubly weighted estimator developed by \citet{4871} for exposure stratified case-cohort studies is intriguing in that it uses 
a separate set of (time-dependent) weights for each covariate.
A preliminary analysis is conducted to estimate quantities that resemble within stratum conditional expectations of partial score contributions given the phase one data, and these are used to form the weights.
An extension of their approach to more general two phase stratified sampling designs would be of considerable interest.

This paper is limited in application to semiparametric models that satisfy the rather stringent assumptions \textbf{A1}-\textbf{A4} of \S 2. 
Even in the case of Cox regression, these have been established only under the ``partly unnecessary" conditions imposed by \citet[\S 25.12.1]{4895}.
His assumption that everyone still ``on-study" is censored at the common time $\tau$ would apply to situations in which time $t$ referred to calendar time, everyone was entered on study at $t=0$ and there was a common closing date at $t=\tau$.
It would not apply, however, if subjects were entered on study at various
calendar times but withdrawn on a common closing date, and $t$ was taken to be ``time-on-study".
Nor would it apply if $t$ was ``age" and subjects both entered and exited the study at various ages. 
We look forward to further work that relaxes these assumptions, in particular to a determination as to whether or not the general approach extends to Cox regression with time-dependent covariates and repeated failure events under standard assumptions \citep{3924}.

In his Appendix \citet{4894} remarks
\begin{quote}
``To our knowledge, there does not exist a general theory on the conditions required for the tightness and weak convergence of Horvitz-Thompson processes. 
However, the results of \citet[\S\S 2.9, 3.6, 3.7]{4920} can be applied to possibly stratified simple random sampling and can potentially be extended to other survey designs."
\end{quote}
One purpose of this paper has been to carry out in detail the program 
mentioned for stratified random sampling.
We conjecture that our fundamental equation (\ref{eq:general}) applies 
to Horvitz-Thompson estimators for other complex sampling designs, and work is in progress to explore these extensions.

\section*{Acknowledgements}
The second author owes thanks to Galen Shorack for a helpful discussion concerning 
the representation in Appendix A.
Supported in part by grants 5-R01-CA40644 and 2-R01-AI291968 from the 
U.S. National Institutes of Health and by grant DMS-0503822 
from the U.S. National Science Foundation.

\bibliographystyle{ims}
\bibliography{sjs}

\section{Appendices}
In Appendices A and B we establish two results slightly more
general than needed for the development in Section 4.   
(See the end of Appendix  B for the special case required.)  
The notation in these two appendices should be understood
to be independent of the that in the body of the paper.

\par\noindent
{\bf Appendix A.  \ A Representation of Stratified Sampling.}
 
 Suppose that $( \Omega, {\cal A},P)$ is a probability space
and $W : (\Omega , {\cal A} ) \rightarrow ({\cal W} , {\cal B} )$. 
Write $P^W$ for the  measure induced by $W$ on 
$ ({\cal W} , {\cal B} )$; in the notation of section 2, $P^W = \tilde{P}_0$.  
Suppose that ${\cal W}_1 , \ldots , {\cal W}_J$ is a (measurable) partition
of ${\cal W}$:  \\
(a) \ \ ${\cal W}_j \in {\cal B}$, $j = 1, \ldots , J$;\\
(b) \ \  ${\cal W}_j \cap {\cal W}_{j'} = \empty$ for $j \not= j'$; and \\ 
(c) \ \ $\cup_{j=1}^J {\cal W}_j = {\cal W}$.\\
We will assume that $P(W \in {\cal W}_j ) \equiv p_j > 0$
for $j=1, \ldots , J$.

Now consider a new probability space $(\Omega^{\dagger}, \cal A^{\dagger}, P^{\dagger})$
where
\begin{eqnarray*} 
&&\Omega^{\dagger}
= \Omega_0^{\dagger} \times \Omega_1^{\dagger} \times \cdots \times \Omega_J^{\dagger} ,\\
&& {\cal A}^{\dagger} 
= {\cal A}_0^{\dagger} \times {\cal A}_1^{\dagger} \times \cdots \times {\cal A}_J^{\dagger} , \\
&& P^{\dagger} = P_0^{\dagger} \cdot P_1^{\dagger} \cdots P_J^{\dagger} , 
\end{eqnarray*}
and random variables $\Delta = (\Delta_1, \ldots , \Delta_J)$, 
$W_1^{\dagger}, \ldots , W_J^{\dagger}$
defined thereon as follows:
for $\omega^{\dagger} = (\omega_0^{\dagger} , \omega_1^{\dagger}, \ldots , 
\omega_J^{\dagger} ) \in \Omega^{\dagger}$,
\begin{eqnarray*}
&& \Delta (\omega^{\dagger}) 
 = \Delta (\omega_0^{\dagger}) \sim \mbox{Multinomial}_J (1, (p_1 , \ldots , p_J ) ) \\
&& W_j^{\dagger} ( \omega^{\dagger} ) = W_j^{\dagger} (\omega_j^{\dagger} ) 
\sim P_j^{\dagger}
\end{eqnarray*}
for $j =1 , \ldots , J$ where  $p_j = P( W \in {\cal W}_j )$, $j =1 , \ldots , J$, and 
$P_j^{\dagger}$ is defined by 
\begin{eqnarray}
P_j ^{\dagger} (W_j \in B) = \frac{P( W \in B\cap {\cal W}_j)}{P(W \in {\cal W}_j )} 
= \frac{P^W (B \cap {\cal W}_j )}{P^W ( {\cal W}_j )},  
\qquad B \in {\cal B} .
\label{DefnOfPSubJDagger}
\end{eqnarray}
Now define a random variable 
$W^{\dagger} :  (\Omega^{\dagger}, \cal A^{\dagger}) \rightarrow  ( {\cal W} , {\cal B} )$ 
by 
\begin{eqnarray*}
W^{\dagger} (\omega^{\dagger} ) 
= \Delta_1 (\omega_0^{\dagger} ) W_1^{\dagger}( \omega_1^{\dagger} )
       + \cdots +  \Delta_J ( \omega_0^{\dagger} ) X_J^{\dagger} ( \omega_J^{\dagger} )  .
\end{eqnarray*}
Note that $\Delta$, $W_1^{\dagger} , \ldots , W_J^{\dagger}$ are independent by 
construction.
\bigskip

\par\noindent
{\bf Proposition A.1} \ \   $W^{\dagger} \stackrel{d}{=} W$ on $({\cal W}, {\cal B} )$.  
That is, $P^{W^{\dagger}} = P^W $ as measures on $({\cal W}, {\cal B})$.
\medskip

\par\noindent
{\bf Proof.} 
First note that 
\begin{eqnarray}
P^{\dagger} (W^{\dagger}\in {\cal W}_j )
& = & P^{\dagger} (W_j^{\dagger} \in {\cal W}_j , \Delta_j = 1) \nonumber  \\
& = & P^{\dagger} (W_j^{\dagger} \in {\cal W}_j ) P^{\dagger} (\Delta_j = 1) 
          %\nonumber \\
 =  1 \cdot p_j = p_j 
\label{ComputationOfProbDaggerOfFallingInjthElementOfPartition}
\end{eqnarray}
using independence of $\Delta$ and $W_j^{\dagger}$, 
the fact that $W_j^{\dagger}$ takes values in ${\cal W}_j$ with $P^{\dagger}$-probability $1$, 
and $P^{\dagger} (\Delta_j = 1)= p_j $ by the definition of $P^{\dagger}$.

Now let $B \in {\cal B}$.  Then since $p_j > 0$ for $j =1, \ldots , J$,
\begin{eqnarray*}
P^{\dagger} (W^{\dagger} \in B) 
& = & \sum_{j=1}^J P^{\dagger} ( W^{\dagger} \in B \cap {\cal W}_j ) 
 =  \sum_{j=1}^J \frac{P^{\dagger} ( W^{\dagger} \in B \cap {\cal W}_j ) }{ P^{\dagger} (W^{\dagger} \in {\cal W}_j )}
                 P^{\dagger} (W^{\dagger} \in {\cal W}_j )\\
& = & \sum_{j=1}^J \frac{P^{\dagger} ( W_j^{\dagger} \in B ) }{P^{\dagger} (W^{\dagger}_j \in {\cal W}_j )} p_j
\qquad \mbox{by} \ (\ref{ComputationOfProbDaggerOfFallingInjthElementOfPartition}) \\
& = & \sum_{j=1}^J \frac{P^W (B \cap {\cal W}_j) / P^W ( {\cal W}_j )}{ 1} \cdot p_j
\qquad \mbox{by} \ (\ref{DefnOfPSubJDagger})  \\
& = & \sum_{j=1}^J P^W ( B \cap {\cal W}_j ) = P^W (B) = P( W \in B) .  
%\qquad \qquad \qquad \qquad\qquad \Box
\end{eqnarray*}
\hfill $\Box$
\medskip

If $W_1, \ldots , W_N$ are i.i.d. $P^W$, then we can represent the $W_i$'s 
in terms of $( \Delta_i , W_{1,i}^{\dagger} , \ldots, W_{J,i}^{\dagger} )$, $i=1, \ldots , N$, i.i.d. as  
$( \Delta , W_1^{\dagger} , \ldots , W_J^{\dagger})$ as described in proposition A.1.  It follows that 
\begin{eqnarray}
\PP_{j, N_j} 
& = & \frac{1}{N_j} \sum_{i=1}^N \delta_{W_i} 1_{{\cal W}_j} (W_i) \nonumber \\ 
& = & \frac{1}{N_j} \sum_{j'=1}^J \sum_{i=1}^N \Delta_{j',i} \delta_{W_{j',i}^{\dagger}} 
            1_{{\cal W}_j} (W_{j,i}^{\dagger}) \nonumber  \\
& = & \frac{1}{N_j} \sum_{i=1}^{N_j} \delta_{W_{j,i} }
\label{eq:StratumSpecificEmpiricalMeasureAppendixB}
\end{eqnarray}
by relabelling the $W_{j,i}^{\dagger}$'s and where 
$N_j = \sum_{i=1}^N \Delta_{j,i}$ on the right side is independent of the $W_{j,i}^{\dagger}$'s.
  This yields the promised doubly indexed form of the stratum - specific
empirical measure in terms of independent $W_{j,i}$'s distributed according to $P_{0|j}$
where, for $B \in {\cal B}$,
$$
P_{0|j} (B) = \frac{P_0 (B 1_{{\cal W}_j} )}{P_0 ( 1_{{\cal W}_j} )} .
$$
\bigskip

\par\noindent
{\bf Appendix B.  Proof of weak convergence of the stratum-specific empirical process}

Let $\PP_{j,N_j} $ be as defined in (\ref{eq:StratumSpecificEmpiricalMeasureAppendixB})
%eq:stratumwiseEmpirical}),
%where the double subscripting:
$$
\PP_{j,N_j} = \frac{1}{N_j} \sum_{i=1}^n \delta_{W_i} 1_{{\cal W}_j} (W_i) 
$$
where 
$$
N^{-1} N_j = \PP_N ( 1_{{\cal W}_j} ) \rightarrow_{a.s.} P_0 ( {\cal W}_j ) \equiv \nu_j > 0 .
$$
\medskip

\par\noindent
{\bf Proposition B.1.}
If ${\cal F}$ is $P_{0}-$Donsker and $\nu_j > 0$, then 
${\cal F}$ is $P_{0|j} - $Donsker on stratum ${\cal W}_j$  in the sense that 
\begin{equation}
\GG_{j,Nj} \equiv \sqrt{N_j} ( \PP_{j, N_j} - P_{0|j} ) 
\rightsquigarrow  \GG_j  \qquad \mbox{in} \ \ \ell^{\infty} ({\cal F}) 
\label{EmpiricalProcessForStratumj}
\end{equation}
where $\GG_j$, defined by 
\begin{equation}
\GG_j (f) = \nu_j^{-1/2} \GG_{P_0} ((f - P_{0|j} (f)) 1_{{\cal W}_j} ) , 
\qquad f \in \ell^{\infty} ({\cal F}) ,
\end{equation}
is a $P_{0|j}$-Brownian bridge process.  
\medskip

\par\noindent
{\bf Remark 1.}
Note that 
\begin{eqnarray*}
Var( \GG_{j} (f) ) 
& = & \nu_j^{-1} P_0 \left [ ( f - P_{0|j} (f))^2 1_{{\cal W}_j} \right ]
 =  Var_j (f) \equiv Var(f(W)| W \in {\cal W}_j ).
\end{eqnarray*}
%\medskip

\par\noindent
{\bf Remark 2.}  The proposition implies that the process 
$\sqrt{N_j} ( \PP_{j, N_j} - P_{0|j} ) $ behaves asymptotically the same as that 
of a sample of fixed size drawn from the conditional distribution $P_{0|j}$.    
\medskip

\par\noindent
{\bf Proof of the proposition}.  First proof.
By the discussion at the beginning of section 2.10.4, page 200, van der Vaart 
and Wellner (1996), ${\cal F}_j \equiv \{ f 1_{{\cal W}_j} : \ f \in {\cal F} \}$ 
is $P_0-$Donsker, and hence the collection 
$\tilde{{\cal F}}_j \equiv \{ f 1_{{\cal W}_j} : \ f \in {\cal F} \cup \{ 1\} \}$
is also $P_0-$Donsker.  Now we 
write 
\begin{eqnarray*}
\sqrt{N_j}  ( \PP_{j,N_j} f - P_{0|j} f ) 
& = & \sqrt{N_j} \left ( \frac{\frac{1}{N} \sum_{i=1}^N f(W_i ) 1_{{\cal W}_j} (W_i)}
                                         {\frac{1}{N} \sum_{i=1}^N  1_{{\cal W}_j} (W_i)}
                                         - \frac{P_0 ( f 1_{{\cal W}_j}}{P_0 ( 1_{{\cal W}_j}) } 
                         \right )  \\
& = & \sqrt{ \frac{N_j}{N}} \left \{  \frac{ \GG_N ( f 1_{{\cal W}_j} ) }{N_j/N}  
           -   \frac{ \GG_N ( 1_{{\cal W}_j} ) P_0 ( f 1_{{\cal W}_j} )}
                                                   { (N_j /N) P_0 ( {\cal W}_j ) }  \right \} \\
& = &  \frac{1}{\sqrt{N_j/N}} \left \{ \GG_N ( f 1_{{\cal W}_j} )  
             -   \GG_N ( 1_{{\cal W}_j} ) P_{0|j} ( f )    \right \} \\
& = & \frac{1}{\sqrt{N_j/N}} \GG_N ( ( f   - P_{0|j} (f)) 1_{{\cal W}_j}  )  \\
& \Rightarrow &   \frac{1}{\sqrt{\nu_j}} \GG_{P_0} ( ( f   - P_{0|j} (f)) 1_{{\cal W}_j}  ) 
 \equiv \GG_{P_{0|j}} (f) \,  ,   
%\qquad \Box                             
\end{eqnarray*}
and, in fact, 
\begin{eqnarray*}
\left \{ \frac{1}{\sqrt{\nu_j}} \GG_{P_0} ((f - P_{0|j} (f))1_{{\cal W}_j} ) : \ f \in {\cal F} \right \} 
\stackrel{d}{=} \{ \GG_{P_{0|j} } (f) : \ f \in {\cal F} \} .
\end{eqnarray*}

Second proof.  By the second representation of the stratum-specific empirical measure 
$\PP_{j,N_j}$ as  $\PP_{j,N_j} = N_j^{-1} \sum_{i=1}^{N_j} \delta_{W_{j,i}}$ where the 
$W_{j,i}$'s are i.i.d. $P_{0|j}$, it follows that the empirical 
process 
$\GG_{j,N_j} = \sqrt{N_j} ( \PP_{j,N_j} - P_{0|j}) $
is just the empirical process of  i.i.d. $W_{j,i}$'s, but with a random sample size 
$N_j$ independent of the $W_{j,i}$'s.    Since $N_j / N \rightarrow \nu_j > 0$, it follows from 
theorem 3.5.1, page 339, van der Vaart and Wellner (1996), that 
$\GG_{j,N_j} \rightsquigarrow \GG_j$ in $\ell^{\infty} ({\cal F})$  where $\GG_j$ is a
$P_{0|j}-$Brownian bridge process as before.
\hfill$\Box$

In the application of the results of Appendices A and B in section 4 we take 
${\cal W}_1 , \ldots , {\cal W}_J$ to be the measurable partition of ${\cal W}$ induced by 
the partition ${\cal V}_1 , \ldots , {\cal V}_J$ of ${\cal V}$ (i.e. ${\cal W}_j = V^{-1} ({\cal V}_j)$
for $j=1, \ldots , J$ where $V(W) \equiv ( \tilde{X} (X), U)$).    Moreover, the Donsker class 
${\cal F}$ in Proposition B.1 is taken to be a Donsker class of functions of $X$ only 
rather than functions of $W = (X,U)$.  This is exactly what is needed for the development
in section 4.
\medskip

{\bf Appendix C. Proof of equation (\ref{eq:TaylorExpansion}).}
Besides the consistency and asymptotic linearity (\ref{eq:InfFuncAlpha}) 
for $\hat\alpha$ assumed in \S 6, we further assume that $0  <  \sigma  \leq \pi_{\alpha}(v)$ as in (\ref{eq:boundedweights}) and that  
\begin{eqnarray}
\Big | \frac{1}{\pi_{\alpha} (v)} - \frac{1}{\pi_{\alpha_0} (v) }
     - \frac{-\dot{\pi}_0^{T}(v)}{\pi_0^2 (v)} (\alpha - \alpha_0) \Big | 
     \le \psi (v) | \alpha - \alpha_0 |^{1+ \zeta }
\label{eq:DerivativeConditionPlus}
\end{eqnarray}
for $\alpha$ in a neighborhood of $\alpha_0$ where $\zeta> 0$ 
and $\psi$ satisfies $E \psi^2 (V) < \infty$.
The second assumption will typically follow from the first 
provided that $\pi_{\alpha}$ has a continuous second derivative.
For example, suppose that $\pi_{\alpha}$ is given by a logistic regression model with linear predictor 
$\tilde v^{T}\alpha$ where $\tilde v=\tilde v(v) \in \RR^q$. 
Then Taylor's formula with remainder shows that the LHS of (\ref{eq:DerivativeConditionPlus}) equals
$\left | \frac{1}{2}e^{-\tilde v^{T} \alpha^*}(\alpha-\alpha_0)^{T}\tilde v \tilde v^{T}(\alpha-\alpha_0) \right |$
with $\alpha^*$ on the line segment between $\alpha$ and $\alpha_0$.
Thus the condition holds with $\zeta=1$ provided $e^{\tilde v^{T} \alpha}=\pi_{\alpha}(v)/[1-\pi_{\alpha}(v)]$ 
is bounded away from 0 and $\tilde V$ has finite fourth moment.
It follows that 
\begin{eqnarray}
\left ( \PP_N^{\hat\pi}- \PP_N^{\pi_0} \right ) \tilde\ell_0 
& = & 
\frac{1}{N}\sum_{i=1}^N {\bf 1}_{{\cal V}_0^c} (V_i) \left ( \frac{\xi_i}{\hat\pi_i} - \frac{\xi_i}{\pi_0} \right ) \tilde\ell_0(X_i) 
         \nonumber \\
& = & \frac{1}{N}\sum_{i=1}^N {\bf 1}_{{\cal V}_0^c}(V_i) 
          \xi_i \tilde\ell_0(X_i) 
           \left [\frac{1}{\pi_{\hat\alpha}(V_i)} -  \frac{1}{\pi_{\alpha_0}(V_i)} 
           - \frac{-\dot{\pi}_0^T(V_i)}{\pi_0^2 (V_i)} (\hat{\alpha} - \alpha_0)
           \right ]   \nonumber \\          
 && \qquad  + \frac{1}{N}\sum_{i=1}^N {\bf 1}_{{\cal V}_0^c}(V_i) 
           \xi_i \tilde\ell_0(X_i) 
           \left [ \frac{-\dot{\pi}_0^T(V_i)}{\pi_0^2 (V_i)}
           \right ]  ( \hat{\alpha} - \alpha_0)  \nonumber \\
& \equiv & R_N -   \frac{1}{N}\sum_{i=1}^N {\bf 1}_{{\cal V}_0^c}(V_i) 
          \frac{\xi_i}{\pi_{0}(V_i)} \tilde\ell_0(X_i) 
           \left [ \frac{\dot{\pi}_0^T(V_i)}{\pi_0 (V_i)} \right ] (\hat{\alpha} - \alpha_0)  
\label{eq:RemainderPlusMainTerm}
\end{eqnarray}
 where by (\ref{eq:boundedweights}), the similar assumption for 
 $\pi_{\alpha}$  and (\ref{eq:DerivativeConditionPlus}),
 \begin{eqnarray*}
 |R_N | 
 & \le & \Big |  \frac{1}{N}\sum_{i=1}^N {\bf 1}_{{\cal V}_0^c}(V_i) 
          \xi_i \tilde\ell_0(X_i) 
           \left [\frac{1}{\pi_{\hat\alpha}(V_i)} -  \frac{1}{\pi_{\alpha_0}(V_i)} 
           - \frac{-\dot{\pi}_0^T(V_i)}{\pi_0^2 (V_i)} (\hat{\alpha} - \alpha_0) \right ]
       \Big | \\
 & \le & \frac{1}{\sigma^2} \frac{1}{N} \sum_{i=1}^N \psi (V_i)|  \tilde\ell_0 (X_i)| \cdot | \hat{\alpha} - \alpha_0|^{1+\zeta} \\
 & = &   O_p (1) | \hat{\alpha} - \alpha_0| | \hat{\alpha} - \alpha_0|^{\zeta} \\
 & = &   O_p (1) O_{p} (N^{-1/2} ) o_p (1).
 \end{eqnarray*}
Multiplying through (\ref{eq:RemainderPlusMainTerm}) by $\sqrt{N}$, we conclude that (\ref{eq:TaylorExpansion}) holds by virtue of $\sqrt{N} \tilde{R}_N = o_p (1)$ and the strong law of large numbers.

\end{document}